\documentclass[11pt]{article}
\usepackage{amsmath,amssymb,amsthm,geometry,verbatim,enumerate,float,tikz}
\newtheorem{thm}{Theorem}
\newtheorem{lemma}[thm]{Lemma}

\usepackage{setspace}

\begin{document}

\onehalfspace

\title{A Characterization of Substar Graphs}

\author{Felix Joos\footnote{Institut f\"{u}r Optimierung und Operations Research, 
Universit\"{a}t Ulm, Ulm, Germany,
e-mail: \texttt{felix.joos@uni-ulm.de}}}

\date{}

\maketitle


\begin{abstract}
The intersection graphs of stars in some tree are known as substar graphs.
In this paper we give a characterization of substar graphs by the list of minimal forbidden induced subgraphs.
This corrects a flaw in the main result of 
Chang, Jacobson, Monma and West (Subtree and substar intersection numbers, Discrete Appl. Math. \textbf{44}, 205-220 (1993))
and this leads to a different list of minimal forbidden induced subgraphs.

\bigskip


\end{abstract}

\section{Introduction}
A graph $G$ is \textit{chordal}, 
if it has no induced cycle of order at least four.
Chordal graphs are one of the most fundamental graph classes and well investigated \cite{bls}.
For a family $\mathcal{M}$ of sets,
an $\mathcal{M}$-\textit{intersection representation} of a graph $G$
is a function $f: V(G)\rightarrow \mathcal{M}$ such that two distinct vertices $x$ and $y$ are adjacent if and only if
$f(x)\cap f(y)\not= \emptyset$.
If $G$ has an $\mathcal{M}$-intersection representation, 
then $G$ is an $\mathcal{M}$-\textit{graph}.

It is a classic result \cite{gavril2} in graph theory
that every chordal graph is a $\mathcal{T}$-graph,
where $\mathcal{T}$ is the set of all subtrees of some tree.
There are several important and well investigated subclasses of chordal graphs.
Let $\mathcal{T}_P$ be the set of all paths in some tree.
We call a $\mathcal{T}_P$-graph a \textit{path graph}.
Path graphs are investigated in many different varieties \cite{gavril,momna}.
If $\mathcal{T}_I$ is the set of all subpaths of some path,
then the class of $\mathcal{T}_I$-graphs are known as \textit{interval graphs},
which are another very important graph class \cite{fishburn}.

A \textit{star} is a tree such that there is a vertex that is adjacent to all other vertices of the tree.
Let $\mathcal{T}_S$ be the set of all substars in some tree $T$.
We call a $\mathcal{T}_S$-graph a \textit{substar graph}.
Chang et al. \cite{cjmw} claimed to characterize substar graphs by a finite list of forbidden induced subgraphs.
Since there are infinitely many minimal forbidden induced subgraphs for substar graphs, 
their claim is false.
In this paper we give a characterization of substar graphs by the infinite list of minimal forbidden induced subgraphs.

Cerioli and Szwarcfiter \cite{cs} characterized a subclass of substar graphs,
which is the class of starlike graphs.
These are intersection graphs of substars of a star.

We start by introducing some basic notation.
Let $G$ be a graph.
We denote by $V(G)$ and $E(G)$ the vertex set and the edge set of $G$, respectively.
A subset $Q$ of the vertices of $G$ is a \textit{clique} in $G$ if every pair of vertices of $Q$ are adjacent, and
$Q$ is a \textit{maximal clique} in $G$ if it is a clique not properly contained in another clique of $G$.
Let the number of edges in a longest shortest path in $G$ be the \textit{diameter} of $G$.
For $k\in \mathbb{N}$, we write $[k]$ for the set $\{1,\ldots,k\}$.

\section{Results}

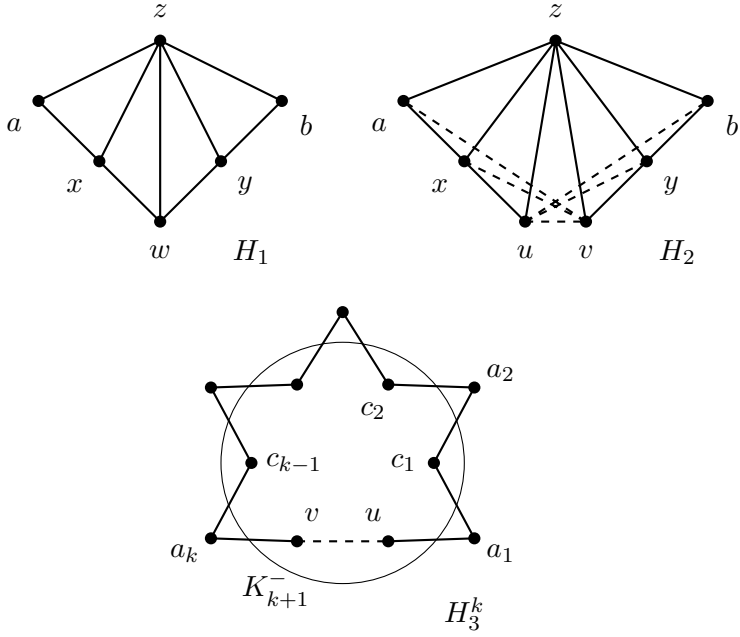
\begin{figure}[t]
\begin{center}
\begin{tikzpicture}[scale=0.8]

\def\ver{0.1} 
\def\x{1}

\def\xa{-3}
\def\ya{4}

\def\xb{3}
\def\yb{4}

\path[fill] (\xa,\ya) circle (\ver);
\path[fill] (\xa-1,\ya+1) circle (\ver);
\path[fill] (\xa-2,\ya+2) circle (\ver);
\path[fill] (\xa+1,\ya+1) circle (\ver);
\path[fill] (\xa+2,\ya+2) circle (\ver);
\path[fill] (\xa,\ya+3) circle (\ver);

\draw[thick] (\xa,\ya)--(\xa-2,\ya+2)
(\xa,\ya)--(\xa+2,\ya+2)
(\xa,\ya)--(\xa,\ya+3)
(\xa-1,\ya+1)--(\xa,\ya+3)
(\xa-2,\ya+2)--(\xa,\ya+3)
(\xa+1,\ya+1)--(\xa,\ya+3)
(\xa+2,\ya+2)--(\xa,\ya+3);

\node (1) at (\xa,\ya-0.5) {$w$};
\node (1) at (\xa-1.4,\ya+0.6) {$x$};
\node (1) at (\xa-2.4,\ya+1.6) {$a$};
\node (1) at (\xa+1.4,\ya+0.6) {$y$};
\node (1) at (\xa+2.4,\ya+1.6) {$b$};
\node (1) at (\xa,\ya+3.5) {$z$};
\node (1) at (\xa+1.5,\ya-0.5) {$H_1$};

\path[fill] (\xb,\yb) circle (\ver);
\path[fill] (\xb-1,\yb+1) circle (\ver);
\path[fill] (\xb-2,\yb+2) circle (\ver);
\path[fill] (\xb+1,\yb) circle (\ver);
\path[fill] (\xb+2,\yb+1) circle (\ver);
\path[fill] (\xb+3,\yb+2) circle (\ver);
\path[fill] (\xb+0.5,\yb+3) circle (\ver);

\draw[thick] (\xb,\yb)--(\xb-2,\yb+2)
(\xb+1,\yb)--(\xb+3,\yb+2)
(\xb,\yb)--(\xb+0.5,\yb+3)
(\xb-1,\yb+1)--(\xb+0.5,\yb+3)
(\xb-2,\yb+2)--(\xb+0.5,\yb+3)
(\xb+1,\yb)--(\xb+0.5,\yb+3)
(\xb+2,\yb+1)--(\xb+0.5,\yb+3)
(\xb+3,\yb+2)--(\xb+0.5,\yb+3);

\draw[thick, dashed](\xb,\yb)--(\xb+1,\yb)
(\xb,\yb)--(\xb+2,\yb+1)
(\xb,\yb)--(\xb+3,\yb+2)
(\xb+1,\yb)--(\xb-1,\yb+1)
(\xb+1,\yb)--(\xb-2,\yb+2);

\node (1) at (\xb,\yb-0.5) {$u$};
\node (1) at (\xb-1.4,\yb+0.6) {$x$};
\node (1) at (\xb-2.4,\yb+1.6) {$a$};
\node (1) at (\xb+1,\yb-0.5) {$v$};
\node (1) at (\xb+2.4,\yb+0.6) {$y$};
\node (1) at (\xb+3.4,\yb+1.6) {$b$};
\node (1) at (\xb+0.5,\yb+3.5) {$z$};
\node (1) at (\xb+2.5,\yb-0.5) {$H_2$};

\def\ra{1.5}
\def\rab{1}
\def\rb{2.5}
\def\rbb{3}

\path[draw] (0,0) circle (2);

\path[fill] (0:\ra) circle (\ver);
\path[fill] (60:\ra) circle (\ver);
\path[fill] (120:\ra) circle (\ver);
\path[fill] (180:\ra) circle (\ver);
\path[fill] (240:\ra) circle (\ver);
\path[fill] (300:\ra) circle (\ver);

\path[fill] (30:\rb) circle (\ver);
\path[fill] (90:\rb) circle (\ver);
\path[fill] (150:\rb) circle (\ver);
\path[fill] (210:\rb) circle (\ver);
\path[fill] (330:\rb) circle (\ver);

\draw[thick] (300:\ra)--(330:\rb)
--(0:\ra)--(30:\rb)
--(60:\ra)--(90:\rb)
--(120:\ra)--(150:\rb)
--(180:\ra)--(210:\rb)
--(240:\ra);

\draw[thick, dashed](240:\ra)--(300:\ra);

\node (1) at (300:\rab) {$u$};
\node (1) at (0:\rab) {$c_1$};
\node (1) at (60:\rab) {$c_2$};
\node (1) at (180:\rab-0.2) {$c_{k-1}$};
\node (1) at (240:\rab) {$v$};

\node (1) at (330:\rbb) {$a_1$};
\node (1) at (30:\rbb) {$a_2$};
\node (1) at (210:\rbb) {$a_k$};

\node (1) at (2,-2.5) {$H_3^k$};
\node (1) at (-1.1,-2.1) {$K_{k+1}^-$};

\end{tikzpicture}
\end{center}
\caption{The set $\mathcal{S}$ of minimal forbidden induced subgraphs for substar graphs.}
\end{figure}

Before we state and prove our main result,
we first introduce the set $\mathcal{S}$ of minimal forbidden induced subgraphs for substar graphs.
See Figure 1 for an illustration.
The graph $H_1$ needs no further explanation.
The graphs of type $H_2$ are as in Figure 1 and the dashed edges indicate the following three possibilities:  
\begin{itemize}
	\item $N_{H_2}(u)\cap\{v,y,b\}=N_{H_2}(v)\cap\{u,x,a\}=\emptyset$,
	\item $N_{H_2}(u)\cap\{v,y,b\}=\{v\}$ and $N_{H_2}(v)\cap\{u,x,a\}=\{u,x,a\}$, or
	\item $N_{H_2}(u)\cap\{v,y,b\}=\{v,y,b\}$ and $N_{H_2}(v)\cap\{u,x,a\}=\{u,x,a\}$.
\end{itemize}
The graph $H_3^k$ for some $k\geq 3$ satisfies
$V(H_3^k)=\{u,v,c_1,\ldots,c_{k-1},a_1,\ldots a_k\}$
and 
$$E(H_3^k)=\left(\binom{\{u,v,c_1,\ldots,c_{k-1}\}}{2}\cup \{ua_1,va_k\}\cup\bigcup_{i=1}^{k-1}\{a_{i}c_{i},c_{i}a_{i+1}\}\right)\setminus \{uv\}.$$
Note that the graphs $H_3^k$ for $k\in \{3,\ldots,6\}$ equal the forbidden subgraphs $H_3$, $H_4$, $H_5$, $H_6$, $H_7$ introduced by Chang et al. and moreover, the graphs $H_3$ and $H_6$ are isomorphic.
The graphs $H_3^k$, for $k$ at least 7, are minimal forbidden induced subgraphs for substar graphs and do not contain any graph as an induced subgraph from the list of Chang et al.

As already mentioned above, a chordal graph $G$ has a $\mathcal{T}$-representation,
that is, there is a tree $T$ and a function $f:V(G)\rightarrow \mathcal{T}$,
where $\mathcal{T}$ is the set of all the subtrees of $T$ and two distinct vertices $v,w$ are adjacent if and only if $f(v)\cap f(w)\not=\emptyset$.
We say the tuple $(T,f)$ is a \textit{tree representation} of $G$.
It is well known and not difficult to see that one can choose $T$ such that $f(V(G))=T$ without increasing the diameter of any subtree in $\mathcal{T}$.
Therefore, from now on we may assume that a tree representation $(T,f)$ of a graph $G$ has the property $f(V(G))=T$.

Let $(T,f)$ be a tree representation of a chordal graph $G$ and $t\in V(T)$.
Obviously, $f^{-1}(t)$ is a clique $Q$ in $G$.
Suppose that $Q$ is not a maximal clique of $G$ and let $Q'$ be a maximal clique of $G$ that contains $Q$.
By the Helly property of subtrees of a tree (see for example \cite{bls}) and the maximality of $Q'$,
there is a vertex $t'\in V(T)$ such that $f^{-1}(t')=V(Q')$.
Let $P:tt_1\ldots t_{k}t'$ be the $t,t'$-path in $T$.
Note that $V(Q)\subseteq f^{-1}(t_i)$ for every $i\in[k]$.
Contracting the edge $tt_1$ to a vertex $v_{tt_1}$ leads to a tree $T'$.
Define $f'$ by ${f'}^{-1}(s)=f^{-1}(s)$ if $s\in V(T)\setminus \{t,t_1\}$ and let $f'^{-1}(v_{tt_1})=f^{-1}(t_1)$.
Trivially $(T',f')$ is a tree representation of $G$.
Repeating this process leads to a tree representation $(\tilde{T},\tilde{f})$ of $G$
such that for every vertex $\tilde{t}\in V(\tilde{T})$ the vertex set $\tilde{f}^{-1}(\tilde{t})$ is a maximal clique in $G$.
If there are two vertices $t$ and $t'$ of $\tilde{T}$ such that $\tilde{f}^{-1}(t)=\tilde{f}^{-1}(t')$,
then contracting all edges on the $t,t'$-path in $\tilde{T}$ 
and repeating this process as often as possible
leads to a tree representation $(\hat{T},\hat{f})$ of $G$ such that
there is a bijection between the maximal cliques of $G$ and the vertices of $\hat{T}$.
Note that this process does not increase the diameter of a subtree in the tree representation induced by a vertex of $G$.
Therefore, from now on we may assume that a tree representation $(T,f)$ of a graph $G$ has the property $f(V(G))=T$ and
there is a bijection between the maximal cliques of $G$ and the vertices of $T$.

\noindent
After these preliminaries, we proceed to our results.

\begin{lemma}\label{minimalforbidden}
The graphs in $\mathcal{S}$ are not substar graphs.
\end{lemma}

\noindent
\textit{Proof:}
For contradiction we consider a substar representation $(T,f)$ for some $G\in \mathcal{S}$
such that $f(v)$ does not contain a $P_4$ as a subgraph for every $v\in V(G)$.

First, let $G=H_1$.
There are four maximal cliques $Q_1= \{z,a,x\}$, $Q_2= \{z,x,w\}$, $Q_3=\{z,w,y\}$, and $Q_4=\{ z,y,b\}$
with corresponding vertices $q_1,q_2,q_3$, and $q_4$ in $T$, respectively.
Since $f(v)$ is connected for every $v\in V(G)$,
$q_1q_2\in E(T)$, $q_2q_3\in E(T)$, and $q_3q_4\in E(T)$.
Thus, $T$ and $f(z)$ is a $P_4$, which is a contradiction.

Now suppose that $G$ is a graph of type $H_2$.
There are maximal cliques
which contain the sets $\{z,a,x\}$, $\{z,x,u\}$, $\{z,v,y\}$, and $\{z,y,b\}$, respectively.
Note that these maximal cliques are distinct, 
because for every pair $Q$ and $Q'$ of these sets there is a vertex of $G$ in $Q$
which is nonadjacent to a vertex in $Q'$.
Suppose $f(z)$ is a star.
If the subtree $f(x)$ contains the center vertex of $f(z)$, then $f(y)$ is disconnected,
otherwise $f(x)$ is disconnected.
Thus $f(z)$ is not a star,
which is a contradiction.

Now we suppose $G=H_3^k$ for some $k\geq 3$.
We denote by $Q_u$ and $Q_v$ the two maximal cliques $\{u,c_1,\ldots ,c_{k-1}\}$ and $\{v,c_1,\ldots ,c_{k-1}\}$ with corresponding vertices $q_u$ and $q_v$ in $T$, respectively.
In addition, there are $k$ maximal cliques 
$Q_1=\{u,a_1,c_1\}$, $Q_i=\{c_{i-1},a_i,c_i\}$ for $i\in \{2,\ldots,k-1\}$, and $Q_k=\{c_{k-1},a_k,v\}$ with corresponding vertices $q_1,\ldots,q_k$  in $T$.
Since $f(v)$ is connected for every $v\in V(G)$,
every vertex $w\in Q_u\cap Q_v$ is contained in every maximal clique that is associated with a vertex on the $q_u,q_v$-path in $T$.
If $k=3$ and $q_uq_v$ is not an edge in $T$,
then $q_uq_2$ and $q_2q_v$ are edges in $T$.
Thus $q_uq_1$ is an edge and hence $f(c_1)$ is not a star, which is a contradiction.
Since no  maximal clique beside $Q_u$ and $Q_v$ contains $\{c_1,c_{k-1}\}$ if $k\geq 4$,
we conclude that $q_uq_v$ is an edge in $T$.
Since $f(u)$ and $f(v)$ is connected, we obtain $q_uq_1,q_vq_k\in E(T)$.
Since $f(c_1)$ contains no $P_4$, we conclude $q_uq_2\in E(T)$.
By repeating the last argument several times, we obtain $q_uq_i\in E(T)$ for every $i\in [k]$, and
hence there is a cycle $q_uq_vq_kq_u$ in $T$.
This is a contradiction to our assumption that $T$ is a tree,
which completes the proof.
$\Box$

\bigskip

\noindent
We immediately proceed to our main result.

\begin{thm}\label{mainthm}
A chordal graph is a substar graph if and only if it has no induced subgraph in $\mathcal{S}$.
\end{thm}

\noindent
\textit{Proof:}
By Lemma \ref{minimalforbidden}, a graph in the set $\mathcal{S}$ is not a substar graph.
We show that every graph that is not a substar graph contains a graph from the set $\mathcal{S}$.
For contradiction, we assume that $G$ is a chordal graph
that is a minimal forbidden induced subgraph for substar graphs and does not belong to $\mathcal{S}$.
Thus $G$ is connected.
Let $(T,f)$ be a tree representation of $G$.
We denote vertices of $T$ by $q$ with an optional additional subscript/superscript.
The corresponding maximal clique in $G$ is denoted by $Q$ with the same optional additional subscript/superscript.

We first show that we may assume that every subtree has diameter at most $3$, i.e. no subtree contains a $P_5$,
otherwise this is a contradiction to our choice of $G$.
To see this, take a vertex $u$ that is contained in only one maximal clique $Q$
(every leaf of $T$ contains such a vertex).
Since all proper induced subgraphs of $G$ are substar graphs, 
$G-u$ is a substar graph.
Therefore, let $(T',f')$ be a tree representation of $G - u$.
Since $Q-u$ is a clique,
its vertices are contained in some maximal clique $Q'$.
Add to $T'$ a vertex $q''$ adjacent to $q'$.
Associate to $q''$ all vertices of $Q$.
All subtrees in $T'$ to which a vertex of $G-u$ is associated are still subtrees
and may contain a $P_4$ but not a $P_5$.

After this preliminary consideration, we proceed in the proof of Theorem \ref{mainthm}.
We choose a tree representation $(T,f)$ of $G$ such that the subtree $f(v)$ of some tree $T$ does not contain a $P_5$ for every $v\in V(G)$
and the number of vertices $v$ such that $f(v)$ has diameter $3$ is minimal.
Note that for every edge $qq'\in E(T)$, by our assumptions on $T$ and the connectedness of $G$, we conclude
$Q\setminus Q'\neq \emptyset,$ and
$Q\cap Q'\neq \emptyset$.
Since $G$ is not a substar graph,
there is a vertex $z$ of $G$ such that
$f(z)$ contains a $P_4:q_1q_2q_3q_4$.
Let $q_2q_3$ be the unique middle edge of all $P_4$s in $f(z)$.
Note that for $q_1$ and $q_4$ there are possibly many choices.
We call a neighbor $q$ of $q_2$ other than $q_3$ in $T$ \textit{good}
if there is a vertex $u$ of $G$ such that $u\in Q\cap Q_2$ but $u\notin Q_3$,
and otherwise \textit{bad}.
Let $k$ be a positive integer.
We say that there is a \textit{construction of order} $k$ (with respect to $q_2$ and $q_3$) in $G$ if the following holds:
\begin{itemize}
	\item there is sequence of distinct vertices $w_0,w_1,\ldots,w_k$ of $G$ such that $w_0=z$,
	\item there are distinct neighbors $q^1,\ldots, q^k$ of $q_2$ other than $q_3$ in $T$
\end{itemize}
such that
\begin{itemize}
	\item $f(z)\cap \{q^1,\ldots,q^k\}=\{q^1\}$,
	\item $f(w_i)\cap \{q^1,\ldots,q^k\}=\{q^i,q^{i+1}\}$ for $i\in[k-1]$,
	\item $f(w_k)\cap \{q^1,\ldots,q^k\}=\{q^k\}$, and
	\item $q_3\in f(w_i)$ if and only if $i\in \{0,\ldots, k-1\}$.
\end{itemize}
We say that there is a \textit{mirrored construction of order} $k$ (with respect to $q_2$ and $q_3$) in $G$, 
if there is a construction of order $k$ with respect to $q_2$ and $q_3$ in $G$
after renaming $q_i$ by $q_{5-i}$ for $i\in[4]$. 
We say that we \textit{switch} a neighbor $q$ of $q_2$ other that $q_3$,
if we delete the edge $qq_2$ and add the edge $qq_3$ in $T$.

In the following part of the proof we show that there is a construction of order $k$
and a mirrored construction of order $\ell$ with respect to $q_2$ and $q_3$ in $G$ for some positive integers $k$ and $\ell$.
Afterwards we show that this implies the existence of an induced subgraph in $G$ from the set $\mathcal{S}$.

If there is a good neighbor of $q_2$ in $f(z)$,
then there is a construction of order $1$ in $G$.
Therefore,
we may assume that all neighbors of $q_2$ in $f(z)$ other that $q_3$ are bad
and denote this set by $V_1$.
Now we simultaneously switch all vertices of $V_1$ and get a new tree representation $(T_1,f_1)$ of $G$.
If no diameter of a subtree $f(v)$ ($v\in V(G)$) increased by this switch, 
then the subtree $f_1(z)$ is a star and the number of vertices $v$ of $G$ such that $f_1(v)$ has diameter $3$ decreased,
which is a contradiction to our choice of $T$.
Let $V_1'\subseteq V_1$ be the set of vertices of $V_1$ 
that contains vertices $w$ of $G$ such that the diameter of $f_1(w)$ is larger than the diameter of $f(w)$.
This implies that $V_1'$ is exactly the set of vertices of $V_1$ in $T_1$,
such that there is for every $q\in V_1'$ a vertex $x_q$ of $G$ with the following properties
\begin{itemize}
	\item $q\in f(x_q)$, 
	\item $q_2,q_3\in f(x_q)$, 
	\item there is a neighbor $\check{q}$ of $q_2$ other than $q_3$ in $T_1$ such that $\check{q}\in f(x_q)$, and
	\item there is no neighbor $q'$ of $q_3$ other than $q_2$ and other than the vertices of $V_1$ in $T_1$ such that $x_q\in Q'$.
\end{itemize}
Note that for every $q\in V_1'$ 
the diameter of $f(x_q)$ increases from 2 to 3 or from 3 to 4 by switching $q$.
%
Let $V_2$ be set of all neighbors $\hat{q}$ of $q_2$ in $T_1$ other that $q_3$ with the property that
there is a vertex $x$ of $G$ such that
\begin{itemize}
	\item $\hat{q}\in f(x)$,
	\item $\tilde{q}\in f(x)$ for some $\tilde{q}\in V_1'$, 
	(note that these properties imply that $q_2,q_3\in f(x)$), and 
	\item there is no neighbor $q'$ of $q_3$ other than $q_2$ and other than the vertices of $V_1$ in $T_1$ such that $q'\in f(x_q)$.
\end{itemize}
That means in $V_2$ are all vertices of $T_1$ that make problems because we switched the vertices of $V_1$.
If there is a good vertex in $V_2$,
then there is construction of order $2$ in $G$.
Therefore,
we may assume all neighbors of $q_2$ in $V_2$ are bad.
Now we simultaneously switch all vertices of $V_2$ and get a new tree representation $(T_2,f_2)$ of $G$.
Note that for every vertex $v$ of $G$ for which the diameter of $f_1(v)$ is larger than the diameter of $f(v)$
the diameter of $f_2(v)$ is as large as $f(v)$ again.
If no diameter of a subtree $f_1(v)$ ($v\in V(G)$) increased by this switch, 
then the subtree $f_2(z)$ is a star and the number of vertices $v$ of $G$ such that $f_2(v)$ has diameter $3$ decreased,
which is a contradiction to our choice of $T$.
Let $V_2'\subseteq V_2$ be the set of vertices of $V_2$ 
that contain vertices $w$ of $G$ such that the diameter of $f_2(w)$ is larger than the diameter of $f_1(w)$.
This implies that $V_2'$ is exactly the set of vertices of $V_2$ in $T_2$,
such that there is for every $q\in V_2'$ a vertex $x_q$ of $G$ with the following properties
\begin{itemize}
	\item $q\in f(x_q)$, 
	\item $q_2,q_3\in f(x_q)$, 
	\item there is a neighbor $\check{q}$ of $q_2$ other than $q_3$ in $T_2$ such that $\check{q}\in f(x_q)$, and
	\item there is no neighbor $q'$ of $q_3$ other than $q_2$ and other than the vertices of $V_2$ in $T_2$ such that $q'\in f(x_q)$.
\end{itemize}
Note that for every $q\in V_2'$ 
the diameter of $f(x_q)$ increases from 2 to 3 or from 3 to 4 by switching $q$.

It is easy to repeat this arguments for $V_i,V_i'$ for $i>2$,
if $V_i,V_i'$ is defined analogously.
Extending the argumentation from above,
we know, if there is a good vertex $v\in V_k$ for some $j\in \mathbb{N}$,
then there is a construction of order $j$ in $G$.
However, if all vertices in $V_j$ are bad, 
then $V_{j+1}$ is nonempty.
Since the number of vertices is finite,
there is a good vertex in $V_k$ for some positive integer $k$.
Thus we obtain a construction of order $k$ in $G$.
By symmetry, we can argue analogously, if we swap the names of $q_2$ and $q_3$ .
Thus we obtain a construction of order $k$ 
and a mirrored construction of order $\ell$ for some positive integers $k,\ell$ with respect to $q_2$ and $q_3$ in $G$.

\noindent
Suppose $k=1$ and $\ell=1$,
then let 
\begin{itemize}
	\item $x$ be the vertex that is associated with $w_1$ in the construction of order $1$,
	\item $y$ be the vertex that is associated with $w_1$ in the mirrored construction of order $1$,
	\item $a\in Q_1\setminus Q_2$,
	\item $b\in Q_4\setminus Q_3$,
	\item $u\in Q_2\setminus Q_1$, and
	\item $v\in Q_3\setminus Q_4$.
\end{itemize}
If $u=v$, then these vertices, together with $z$, induce the graph $H_1$,
otherwise the graph $H_1$ or a graph of type $H_2$.

\noindent
Suppose $k+\ell>2$,
then $H_3^{k+\ell}$ is an induced subgraph of $G$.
To see this, let
\begin{itemize}
	\item $z=c_\ell$,
	\item $c_{\ell+i}$ be the vertex that is associated with $w_i$ in the construction of order $k$ for $i\in[k-1]$,
	\item $v$ be the vertex that is associated with $w_k$ in the construction of order $k$,
	\item $c_j$ be the vertex that is associated with $w_{\ell-j}$ in the mirrored construction of order $\ell$ for $j\in[\ell-1]$,
	\item	$u$ be the vertex that is associated with $w_\ell$ in the mirrored construction of order $\ell$,
	\item $a_{\ell+i}$ be a vertex in $Q^i\setminus Q_2$ for $i\in[k]$ where $Q^i$ is associated with the construction of order $k$, and
	\item $a_{\ell+1-j}$ be a vertex in $Q^j\setminus Q_3$ for $j\in[\ell]$ where $Q^j$ is associated with the mirrored construction of order $\ell$.
\end{itemize}
This completes the proof of Theorem \ref{mainthm}.
$\Box$

\end{document}